*Original Article*

# Application of Tikhonov Regularization in Generalized Inverse of Adjacency Matrix of Undirected Graph

Paul Ryan A. Longhas[1] and Alsafat M. Abdul[2]

[1,2]Instructor, Department of Mathematics and Statistics, Polytechnic University of the Philippines, Manila, Philippines

**Abstract** — In this article, we found the Moore-Penrose generalized inverse of adjacency matrix of an undirected graph, explicitly. We proved that the matrix $R_\lambda = [r_{ij}]$ is nonsingular where $r_{ii} = \frac{1}{\lambda} + \deg v_i$ and $r_{ij} = |N_G(v_i) \cap N_G(v_j)|$ for $i \neq j$ and, we proved that $A_G^\dagger = [s_{ij}]_{1 \leq i,j \leq n}$ where $s_{ij} = s_{ji} = \lim_{\lambda \to +\infty} \langle R_\lambda^{-1} e_j, f_i \rangle$. The proof of the main result was based on the Tikhonov regularization.

## I. INTRODUCTION

Let $G = (V, E)$ be an undirected graph of order $n$ and size $m$ where set $V = \{v_1, v_2, \ldots, v_n\}$ is a vertex set, and $E$ is an edge set of $G$. The adjacency matrix of $G$ is an $n \times n$ matrix $A_G = [a_{ij}]_{1 \leq i,j \leq n}$ where

$$a_{ij} = \begin{cases} 1 & \text{if } v_i v_j \in E \\ 0 & \text{otherwise} \end{cases}. \quad (1)$$

The adjacency matrix $A_G$ is a symmetric matrix, that is, $A_G = A_G^t$ where $A_G^t = [a_{ji}]_{1 \leq i,j \leq n}$ since $a_{ij} = a_{ji}$ [4,5,6,7,8]. The Moore-Penrose generalized inverse of $A_G$ is an $n \times n$ matrix $A_G^\dagger$ where for each $y \in \mathbb{R}^n$ we have the following:
1. $\|A_G A_G^\dagger y - y\| \leq \|A_G x - y\|$ for all $x \in \mathbb{R}^n$.
2. If $\|A_G z - y\| \leq \|A_G x - y\|$ for all $x \in \mathbb{R}^n$, then $\|A_G^\dagger y\| \leq \|z\|$.[1,2,3,24,25]

Here, the norm $\|\cdot\|$ is the standard norm in $\mathbb{R}^n$. Another way to define the Moore-Penrose generalized inverse is using the Moore-Penrose equation, see [2,3,9,10,11,12,24,25].

One of the known approaches in finding the Moore-Penrose generalized inverse of any $m \times n$ matrix $A$ with real number entries is decomposing it using Singular-Value decomposition, say $A = U\Sigma V^t$ where $U$ is an $m \times m$ orthogonal matrix, $\Sigma$ is an $m \times n$ rectangular diagonal matrix with nonnegative real numbers on the diagonal, and $V$ is an $n \times n$ orthogonal matrix. If $A = U\Sigma V^t$ is the Singular-Value decomposition of $A$, then $A^\dagger = V\Sigma^\dagger U^t$ where $\Sigma^\dagger$ is the Moore-Penrose generalized inverse of $\Sigma$, which was formed by replacing every nonzero diagonal entry of $\Sigma$ by its reciprocal and transposing the resulting matrix.[3,18,19]

Another way of computing the Moore-Penrose generalized inverse of any bounded linear transformation was based on the Tikhonov regularization. Given a linear transformation $T: \mathbb{R}^n \to \mathbb{R}^m$, then for $y \in \text{Ran}(T) + \text{Ran}(T)^\perp$, the problem

$$\min_{x \in X}(\alpha\|x\|^2 + \|Tx - y\|^2) \quad (2)$$

has a unique solution $x_\alpha$. Furthermore, $x_\alpha \to T^\dagger y$ as $\alpha \to 0^+$.[1,13,14,15,16,17]

In this study, we found the Moore-Penrose generalized inverse of the adjacency matrix of the undirected graph $G$ using the Tikhonov regularization. Specifically, we proved that if $A_G$ is the adjacency matrix of $G$, $\lambda > 0$ and $R_\lambda = [r_{ij}]_{1 \leq i,j \leq n}$ where

$$r_{ij} = \begin{cases} |N_G(v_i) \cap N_G(v_j)| & \text{if } i \neq j \\ \dfrac{1}{\lambda} + \deg v_i & \text{if } i = j \end{cases} \quad (3)$$

then $R_\lambda$ is nonsingular, for all $\lambda > 0$ and $A_G^\dagger = [s_{ij}]_{1 \leq i,j \leq n}$ where

$$s_{ij} = s_{ji} = \lim_{\lambda \to +\infty} \langle R_\lambda^{-1} e_j, f_i \rangle. \quad (4)$$

In this paper, $G$ is always an undirected graph, the norm $\|\cdot\|$ is the standard norm in $\mathbb{R}^n$, the inner product $\langle \cdot, \cdot \rangle$ is the standard inner product in $\mathbb{R}^n$, and the set $\{e_1, e_2, \dots, e_n\}$ is the standard basis of $\mathbb{R}^n$. In addition to that, the set of neigbourhood of $v \in V$ is denoted by $N_G(v) = \{u : uv \in E\}$.

## II. PRELIMINARY

Let $G = (V, E)$ be an undirected graph of order $n$ with size $m$ where $V = \{v_1, v_2, v_3, \dots, v_n\}$ and $E$ be the edge set of $G$. For $i = 1,2,3,\dots,n$, set

$$\delta_u(v) = \begin{cases} 1 & \text{if } uv \in E \\ 0 & \text{if } uv \notin E \end{cases} \quad (5)$$

and vector $f_i = \left(\delta_{v_i}(v_1), \delta_{v_i}(v_2), \delta_{v_i}(v_3), \dots, \delta_{v_i}(v_n)\right) \in \mathbb{R}^n$. Then, the following proposition is easy to stablish:

**Proposition 2.1.** The following hold:
i. If $A_G$ is the adjacency matrix of $G$, then for every $x \in \mathbb{R}^n$, we have

$$A_G x = \sum_{i=1}^n \langle x, f_i \rangle e_i. \quad (6)$$

ii. For every $y \in \mathbb{R}^n$, we have

$$A_G^t y = \sum_{i=1}^n \langle y, e_i \rangle f_i. \quad (7)$$

iii. For $i = 1,2,3,\dots n$, we have $A_G e_i = f_i$.
iv. For all $i, j \in \{1,2,3,\dots,n\}$, we have $\langle f_i, f_j \rangle = |N_G(v_i) \cap N_G(v_j)|$. Furthermore, for $i = 1,2,3,\dots,n$ we have $\|f_i\|^2 = \deg v_i$.

Proof: The proof for (i) and (iv) are straightforward. Thus, we will only prove (ii) and (iii).
Proof of (ii): Let $x, y \in \mathbb{R}^n$. Then,

$$\langle A_G^t y, x \rangle = \langle y, A_G x \rangle = \langle y, \sum_{i=1}^n \langle x, f_i \rangle e_i \rangle = \sum_{i=1}^n \langle x, f_i \rangle \langle y, e_i \rangle = \langle \sum_{i=1}^n \langle y, e_i \rangle f_i, x \rangle$$

Therefore,

$$A_G^t y = \sum_{i=1}^n \langle y, e_i \rangle f_i. \quad (8)$$

Proof (iii): Let $i \in \{1,2,3,\dots n\}$. Since $A_G = A_G^t$, then

$$A_G e_i = A_G^t e_i = \sum_{k=1}^n \langle e_i, e_k \rangle f_k = f_i$$

since $\langle e_i, e_j \rangle = \delta_{ij}$ where $\delta_{ij}$ is the Kronecker delta. □

## III. MAIN RESULTS

It was well-known from the Tikhonov regularization that if $A$ is an $n \times n$ matrix with real entries, then for $y \in \mathbb{R}^n$, the problem

$$\min_{x \in \mathbb{R}^n} (\alpha \|x\|^2 + \|Ax - y\|^2) \quad (9)$$

has unique solution $x_\alpha$. Furthermore, $\|x_\alpha - A^\dagger y\| \to 0^+$ [1,13,14,15,16,17]. The next theorem was based on the Tikhonov regularization.

**Theorem 3.1** Let $G = (V, E)$ be undirected graph of order $n$ and size $m$. For $\lambda > 0$, define the $n \times n$ matrix $R_\lambda = [r_{ij}]_{1 \leq i,j \leq n}$ where

$$r_{ij} = \begin{cases} |N_G(v_i) \cap N_G(v_j)| & \text{if } i \neq j \\ \frac{1}{\lambda} + \deg v_i & \text{if } i = j \end{cases} \quad (10)$$

Then, $R_\lambda$ is nonsingular, for all $\lambda > 0$ and $A_G^\dagger = [s_{ij}]_{1 \leq i,j \leq n}$ where $s_{ij} = s_{ji} = \lim \langle R_\lambda^{-1} e_j, f_i \rangle$.

**Proof:** Observe that $\frac{1}{\lambda} I_n + A_G^t A_G = \frac{1}{\lambda} I_n + A_G^2 = R_\lambda$. Since $A_G^2$ is a positive operator, then $-\frac{1}{\lambda}$ is not an eigenvalue of $A_G^2$ for all $\lambda > 0$, and thus, $R_\lambda = \frac{1}{\lambda} I_n + A_G^2$ is nonsingular. Let $y \in \mathbb{R}^n$, and let $x_\lambda$ be the unique solution to the problem:

$$\min_{x \in \mathbb{R}^n} \left( \frac{1}{\lambda} \|x\|^2 + \|A_G x - y\|^2 \right) \quad (11)$$

Then, we have

$$\nabla \left( \frac{1}{\lambda} \|x_\lambda\|^2 + \|A_G x_\lambda - y\|^2 \right) = 0 \quad (12)$$

The equation in (12) implies that $\left( \frac{1}{\lambda} I_n + A_G^2 \right) x_\lambda = A_G y$, and therefore, $R_\lambda x_\lambda = A_G y$. Since $R_\lambda$ is nonsingular, then $x_\lambda = R_\lambda^{-1} A_G y$. Tikhonov regularization implies that $R_\lambda^{-1} A_G \to A_G^\dagger$ as $\lambda \to +\infty$, and hence, the matrix representation of $A_G^\dagger$ with respect to the standard basis $\{e_1, e_2, \ldots, e_n\}$ is given by $A_G^\dagger = [s_{ij}]$ where $s_{ij} = \lim_{\lambda \to +\infty} \langle R_\lambda^{-1} e_j, A_G e_i \rangle = \lim_{\lambda \to +\infty} \langle R_\lambda^{-1} e_j, f_i \rangle$. Since $A_G$ is symmetric, then $(A_G^\dagger)^t = (A_G^t)^\dagger = (A_G^\dagger)$. Hence, $s_{ij} = s_{ji}$. □

The next corollary immediately holds.

**Corollary 3.2.** $A_G$ is nonsingular if and only if $\lim_{\lambda \to +\infty} \langle R_\lambda^{-1} f_j, f_i \rangle = \delta_{ij}$, where $\delta_{ij}$ is the Kronecker delta.

**Proof:** $A_G$ is nonsingular if and only if $R_\lambda^{-1} A_G^2 \to I_n$ as $\lambda \to +\infty$ if and only if $\delta_{ij} = \lim_{\lambda \to +\infty} \langle R_\lambda^{-1} A_G e_j, A_G e_i \rangle = \lim_{\lambda \to +\infty} \langle R_\lambda^{-1} f_j, f_i \rangle$, as desired. □

The following are the applications of theorem 3.1 and corollary 3.2.

**Example 1:** *Generalized Inverse of Star Graph*

Let $S_n = (V, E)$ be a star graph of order $n$, where $V = \{v_1, v_2, \ldots, v_n\}$ and $E = \{v_1 v_i : i = 2,3,4, \ldots, n\}$ [20,21]. Then, $A_{S_n}^\dagger = \frac{1}{n-1} A_{S_n}$.

**Proof:** Observe $N_{S_n}(v_i) = \{v_1\}$ for $i = 2,3,4, \ldots, n$, and thus, $|N_{S_n}(v_i) \cap N_{S_n}(v_j)| = 1$ for $i, j \in \{2,3,4, \ldots, n\}$, $\deg(v_i) = 1$ for $i = 2,3,4, \ldots, n$, $\deg(v_1) = n - 1$, and $|N_{S_n}(v_1) \cap N_{S_n}(v_j)| = 0$ for $j = 2,3,4, \ldots, n$. Therefore, we have

$$R_\lambda = \begin{bmatrix} \frac{1}{\lambda} + n - 1 & 0 & 0 & 0 & 0 & & 0 \\ 0 & \alpha & 1 & 1 & 1 & \cdots & 1 \\ 0 & 1 & \alpha & 1 & 1 & & 1 \\ 0 & 1 & 1 & \alpha & 1 & & 1 \\ 0 & & & & & & \\ \vdots & & & & & \ddots & \vdots \\ 0 & 1 & 1 & 1 & 1 & \cdots & \alpha \end{bmatrix}$$

where $\alpha = \frac{1}{\lambda} + 1$. Observe that

$$R_\lambda^{-1} = \begin{bmatrix} \frac{1}{\frac{1}{\lambda} + n - 1} & 0 & 0 & 0 & 0 & \cdots & 0 \\ 0 & \beta & \theta & \theta & \theta & \cdots & \theta \\ 0 & \theta & \beta & \theta & \theta & \cdots & \theta \\ 0 & \theta & \theta & \beta & \theta & \cdots & \theta \\ 0 & \theta & \theta & \theta & \beta & \cdots & \theta \\ \vdots & \vdots & \vdots & \vdots & \vdots & \ddots & \vdots \\ 0 & \theta & \theta & \theta & \theta & \cdots & \beta \end{bmatrix}$$

where $\theta = \frac{1}{-\alpha-n\alpha+3\alpha+n-2}$ and $\beta = \frac{-\alpha-n+3}{-\alpha^2-n\alpha+3\alpha+n-2}$. Thus, $R_\lambda^{-1}e_1 = \left(\frac{1}{\frac{1}{\lambda}+n-1}, 0, 0, 0, \ldots, 0, 0\right)$, and for $k > 1$ we have $R_\lambda^{-1}e_k = (0, \theta, \ldots, \theta, \beta, \theta, \ldots, \theta)$ where the first coordinate is 0, the $k$th coordinate is $\beta$, otherwise $\theta$. Since $f_1 = (0, 1, 1, 1, \ldots, 1)$ and $f_k = (1, 0, 0, \ldots, 0)$ for $k > 1$, then $A_{S_n}^\dagger = [s_{ij}]_{1 \leq i,j \leq n}$ where

$$s_{11} = \lim_{\lambda \to +\infty} \langle R_\lambda^{-1}e_1, f_1 \rangle = \lim_{\lambda \to +\infty} 0 = 0$$

$$s_{1k} = s_{k1} = \lim_{\lambda \to +\infty} \langle R_\lambda^{-1}e_1, f_k \rangle = \lim_{\lambda \to +\infty} \frac{1}{\frac{1}{\lambda}+n-1} = \frac{1}{n-1} \quad (k > 1)$$

$$s_{ij} = s_{ji} = \lim_{\lambda \to +\infty} \langle R_\lambda^{-1}e_i, f_j \rangle = \lim_{\lambda \to +\infty} 0 = 0 \quad (i, j > 1).$$

Therefore,

$$A_{S_n}^\dagger = \begin{bmatrix} 0 & \frac{1}{n-1} & \frac{1}{n-1} & \frac{1}{n-1} & \frac{1}{n-1} & \cdots & \frac{1}{n-1} \\ \frac{1}{n-1} & 0 & 0 & 0 & 0 & \cdots & 0 \\ \frac{1}{n-1} & 0 & 0 & 0 & 0 & \cdots & 0 \\ \frac{1}{n-1} & 0 & 0 & 0 & 0 & \cdots & 0 \\ \frac{1}{n-1} & 0 & 0 & 0 & 0 & \cdots & 0 \\ \vdots & \vdots & \vdots & \vdots & \vdots & \ddots & \vdots \\ \frac{1}{n-1} & 0 & 0 & 0 & 0 & \cdots & 0 \end{bmatrix}$$

**Example 2:** The complete graph of order 4 is the graph $K_4 = (V, E)$, where $V = \{v_1, v_2, v_3, v_4\}$ and $E = \{v_iv_j : i, j \in \{1, 2, 3, 4\}\}$ [22,23]. Since $N_{K_4}(v_k) = \{v_i : i \neq k\}$ for all $k \in \{1, 2, 3, 4\}$, then $\deg(v_k) = 3$ and $|N_{K_4}(v_i) \cap N_{K_4}(v_j)| = 2$ for $i \neq j$. Therefore,

$$R_\lambda = \begin{bmatrix} 3+\frac{1}{\lambda} & 2 & 2 & 2 \\ 2 & 3+\frac{1}{\lambda} & 2 & 2 \\ 2 & 2 & 3+\frac{1}{\lambda} & 2 \\ 2 & 2 & 2 & 3+\frac{1}{\lambda} \end{bmatrix}$$

Taking the inverse of $R_\lambda$, we have

$$R_\lambda^{-1} = \begin{bmatrix} \beta & \theta & \theta & \theta \\ \theta & \beta & \theta & \theta \\ \theta & \theta & \beta & \theta \\ \theta & \theta & \theta & \beta \end{bmatrix}$$

where $\beta = \frac{7+\frac{1}{\lambda}}{(3+\frac{1}{\lambda})(7+\frac{1}{\lambda})-12}$ and $\theta = \frac{-2}{(3+\frac{1}{\lambda})(7+\frac{1}{\lambda})-12}$. Thus, $R_\lambda^{-1}e_1 = (\beta, \theta, \theta, \theta)$, $R_\lambda^{-1}e_2 = (\theta, \beta, \theta, \theta)$, $R_\lambda^{-1}e_3 = (\theta, \theta, \beta, \theta)$ and $R_\lambda^{-1}e_4 = (\theta, \theta, \theta, \beta)$. Since $f_1 = (0, 1, 1, 1), f_2 = (1, 0, 1, 1), f_3 = (1, 1, 0, 1), f_4 = (1, 1, 1, 0)$, then

$$\langle R_\lambda^{-1}e_i, f_i \rangle = 3\theta \quad (i = 1, 2, 3, 4) \tag{13}$$
$$\langle R_\lambda^{-1}e_i, f_j \rangle = \beta + 2\theta \quad (i \neq j). \tag{14}$$

Letting $\lambda \to +\infty$, then $\theta \to \frac{-2}{9}$ and $\beta \to \frac{7}{9}$. Thus, $A_{K_4}^\dagger = [s_{ij}]$ where $s_{ii} = 3\left(-\frac{2}{9}\right) = -\frac{2}{3}$ for $i = 1, 2, 3, 4$ and $s_{ij} = s_{ji} = \frac{7}{9} + 2\left(-\frac{2}{9}\right) = \frac{1}{3}$ for $i \neq j$. Therefore,

$$A_{K_4}^\dagger = \begin{bmatrix} \frac{-2}{3} & \frac{1}{3} & \frac{1}{3} & \frac{1}{3} \\ \frac{1}{3} & \frac{-2}{3} & \frac{1}{3} & \frac{1}{3} \\ \frac{1}{3} & \frac{1}{3} & \frac{-2}{3} & \frac{1}{3} \\ \frac{1}{3} & \frac{1}{3} & \frac{1}{3} & \frac{-2}{3} \end{bmatrix}$$

**Example 3:** Show using corollary 3.2 that the adjacency matrix of $K_4$ is nonsingular, and thus, $A_{K_4}^\dagger = A_{K_4}^{-1}$.

**Proof:** In view of corollary 3.2, we need to show that

$$\lim_{\lambda \to +\infty} \langle R_\lambda^{-1} f_j, f_i \rangle = \delta_{ij} \tag{15}$$

where $\delta_{ij}$ is the Kronecker delta. Applying example 2, note that we have $R_\lambda^{-1} f_1 = (3\theta, \beta + 2\theta, \beta + 2\theta, \beta + 2\theta)$, $R_\lambda^{-1} f_2 = (\beta + 2\theta, 3\theta, \beta + 2\theta, \beta + 2\theta)$, $R_\lambda^{-1} f_3 = (\beta + 2\theta, \beta + 2\theta, 3\theta, \beta + 2\theta)$ and $R_\lambda^{-1} f_4 = (\beta + 2\theta, \beta + 2\theta, \beta + 2\theta, 3\theta)$. Thus,

$$\langle R_\lambda^{-1} f_j, f_i \rangle = \begin{cases} 3(\beta + 2\theta) & \text{if } i = j \\ 3\theta + 2(\beta + 2\theta) & \text{if } i \neq j \end{cases}. \tag{16}$$

Letting $\lambda \to +\infty$, then $3(\beta + 2\theta) \to 1$ and $3\theta + 2(\beta + 2\theta) \to 0$, and therefore, $\lim_{\lambda \to +\infty} \langle R_\lambda^{-1} f_j, f_i \rangle = \delta_{ij}$, as desired. □

## IV. CONCLUSIONS

If $A_G$ is nonsingular, then $A_G^\dagger = A_G^{-1}$. Thus, theorem 3.1 is much difficult to use than by solving the inverse of $A_G$. However, it is useful when the adjacency matrix of the graph is singular since the inverse of $A_G$ does not exists. Consequently, theorem 3.1 states that we can solve the problem on finding the Moore-Penrose generalized inverse of the adjacency matrix of a graph by finding the inverse of the resolvent $R_\lambda = \frac{1}{\lambda} I + A_G^2$ of $A_G^2$, which is always nonsingular for $\lambda > 0$.

## ACKNOWLEDGMENT

The author is grateful to the Department of Mathematics and Statistics of Polytechnic University of the Philippines, Manila for their unending support to finish this paper.